\documentclass[journal]{IEEEtran}
\usepackage{amsmath}
\usepackage{amsfonts}
\usepackage{amssymb}
\usepackage{bbm}
\usepackage{amsthm}
\usepackage{mathtools}
\usepackage{float}
\usepackage{setspace}
\usepackage{graphicx}
\usepackage{pdfpages}
\usepackage{tikz}
\usetikzlibrary{arrows,positioning,decorations.pathreplacing,shapes,arrows,positioning}
\newtheorem{theorem}{Theorem}
\newtheorem{lemma}[theorem]{Lemma}
\newtheorem{problem}[]{Problem}
\newtheorem{proposition}[theorem]{Proposition}

\newtheoremstyle{example}
{3pt}
{3pt}
{}
{}
{\upshape \bfseries}
{.}
{.5em}
{}
\theoremstyle{example}
\newtheorem{example}{Example}
\newtheorem{dfn}{Definition}
\newenvironment{myproof}[1][\proofname]{\proof[\bf #1]}{\endproof}
\def\*#1{\mathbf{#1}}
\tikzstyle{bball} = [circle,shading=ball, ball color=black!100!white,
    minimum size=0.4cm]
\tikzstyle{wball} = [circle,shading=ball, ball color=white!100!black,
    minimum size=0.4cm]
\DeclareMathOperator{\rank}{rank}
\DeclareMathOperator{\im}{im}
\begin{document} 
\title{A Distance-Based Approach to Strong Target Control of Dynamical Networks}
\author{H.~J.~van Waarde, M.~K.~Camlibel and H.~L.~Trentelman%
\thanks{The authors are with the Johann Bernoulli Institute for Mathematics and Computer Science, University of Groningen, P. O Box 800, 9700 AV Groningen, The Netherlands, Email: h.j.van.waarde@student.rug.nl, h.l.trentelman@rug.nl, and m.k.camlibel@rug.nl.}}%


\maketitle

\begin{abstract}
This paper deals with controllability of dynamical networks. It is often unfeasible or unnecessary to fully control large-scale networks, which motivates the control of a prescribed subset of agents of the network. This specific form of output control is known under the name target control. We consider target control of a family of linear control systems associated with a network, and provide both a necessary and a sufficient topological condition under which the network is strongly targeted controllable. Furthermore, a leader selection algorithm is presented to compute leader sets achieving target control. 
\end{abstract}


\IEEEpeerreviewmaketitle

\section{Introduction}
\IEEEPARstart{D}{uring} the last two decades, networks of dynamical agents have been extensively studied. It is customary to represent the infrastructure of such networks by a graph, where nodes are identified with agents and arcs correspond to the communication between agents. In the study of controllability of networks, two types of nodes are distinguished: leaders, which are influenced by external input, and followers whose dynamics are completely determined by the behaviour of their neighbours. Network controllability comprises the ability to drive the states of all nodes of the network to any desired state, by applying appropriate input to the leaders.

Motivated by model uncertainties, the notion of structural controllability of linear control systems fully described by the pair $(A,B)$ was introduced by Lin \cite{lincontrollability}. Here the matrices $A$ and $B$ are zero-nonzero patterns, i.e. each entry of $A$ and $B$ is either a fixed zero or a free nonzero parameter. In this framework, weak structural controllability requires almost all realizations of $(A,B)$ to be controllable. That is: for almost all parameter settings of the entries of $A$ and $B$, the pair $(A,B)$ is controllable. Lin provided a graph-theoretic condition under which $(A,B)$ is weakly structurally controllable in the single-input case. Many papers followed \cite{lincontrollability}, amongst others we name \cite{glover1976} and \cite{multiinput} in which extensions to multiple leaders are given, and the article \cite{strongstructuralcontrollability1979}, that introduces strong structural controllability, which requires all realizations of $(A,B)$ to be controllable. 

In recent years, structural controllability gained much attention in the study of networks of dynamical agents \cite{burgarth2013}, \cite{chapman}, \cite{barabasicontrol}, \cite{zeroforce}, \cite{trefois}. With a given network graph, a family of linear control systems is associated, where the structure of the state matrix of each system depends on the network topology, and the input matrix is determined by the leader set. In this framework, a network is said to be weakly (strongly) structurally controllable if almost all (all) systems associated to the network are controllable. The graph-theoretic results obtained in classical papers \cite{lincontrollability}, \cite{strongstructuralcontrollability1979} lend themselves excellently in the study of structural controllability of networks. A topological condition for weak structural controllability of networks is given in terms of maximum matchings in \cite{barabasicontrol}, while strong structural controllability is fully characterized in terms of zero forcing sets in \cite{zeroforce}. Such graph-theoretic conditions have a considerable advantage in their computational robustness compared to rank conditions, and aid in finding leader selection procedures. 

However, in large-scale networks with high vertex degrees, a substantial amount of nodes must be chosen as leader to achieve full control in the strong sense, which is often unfeasible. Furthermore, in some applications full control over the network is unnecessary. Hence, we are interested in controlling a subset of agents, called target nodes. This specific form of output control is known under the name target control \cite{barabasitarget}, \cite{Nima}. Potential applications of target control within the areas of biology, chemical engineering and economic networks are identified in \cite{barabasitarget}.   

A network is said to be strongly targeted controllable if all systems in the family associated to the network graph are targeted controllable. In this paper we consider strong targeted controllability for the class of state matrices called distance-information preserving matrices. The adjacency matrix and symmetric, indegree and outdegree Laplacian matrices are examples of distance-information preserving matrices. As these matrices are often used to describe network dynamics (see e.g. \cite{egerstedt2012}, \cite{quantumcontrol}, \cite{rahmani2009}, \cite{tanner2004}, \cite{zhang2011}), distance-information preserving matrices form an important class of matrices associated to network graphs.

Our main results are threefold. Firstly, we provide a sufficient topological condition for strong targeted controllability of networks, that substantially generalizes the results of \cite{Nima} for the class of distance-information preserving matrices. Furthermore, we note that the `k-walk theory' described in \cite{barabasitarget} is easily obtained as a special case of our result. Secondly, noting that our proposed sufficient condition for target control is not a one-to-one correspondence, we establish a necessary graph-theoretic condition for strong targeted controllability. Finally, we provide a two-phase leader selection algorithm consisting of a binary linear programming phase and a greedy approach to obtain leader sets achieving target control. 

The paper is organized as follows: in Section \ref{Preliminaries} we introduce preliminaries and notation. Subsequently, the problem is stated in Section \ref{Problem Statement}. Our main results are presented in Section \ref{Main results}. Finally, Section \ref{Conclusions} contains our conclusions.

\section{Preliminaries}
\label{Preliminaries}
Consider a directed graph $G = (V,E)$, where $V$ is a set of $n$ vertices, and  $E$ is the set of directed arcs. The cardinality of a vertex set $V'$ is denoted by $|V'|$. Throughout this paper, all graphs are assumed to be simple and without self-loops. 

We define the distance $d(u,v)$ between two vertices $u,v \in V$ as the length of the shortest path from $u$ to $v$. If there does not exist a path in the graph $G$ from vertex $u$ to $v$, the distance $d(u,v)$ is defined as infinite. Moreover, the distance from a vertex to itself is equal to zero.  

For a nonempty subset $S \subseteq V$ and a vertex $j \in V$, the distance from $S$ to $j$ is defined as
\begin{equation}
d(S,j) := \min_{i \in S} d(i,j).
\end{equation}

A directed graph $G = (V,E)$ is called bipartite if there exist disjoint sets of vertices $V^-$ and $V^+$ such that $V = V^- \cup V^+$ and $(u,v) \in E$ only if $u \in V^-$ and $v \in V^+$. We denote bipartite graphs by $G = (V^-, V^+, E)$, to indicate the partition of the vertex set.

\subsection{Qualitative class and pattern class}
The qualitative class of a directed graph $G$ is a family of matrices associated to the graph. Each of the matrices of this class contains a nonzero element in position $i,j$ if and only if there is an arc $(j,i)$ in $G$, for $i \neq j$. More explicitly, the qualitative class $\mathcal{Q}(G)$ of a graph $G$ is given by
\begin{equation*}
\mathcal{Q}(G) = \{X \in \mathbb{R}^{n \times n} \: | \: \text{ for } i \neq j, \: X_{ij} \neq 0 \iff (j,i) \in E \}.
\end{equation*}

Note that the diagonal elements of a matrix $X \in \mathcal{Q}(G)$ do not depend on the structure of $G$, these are `free elements' in the sense that they can be either zero or nonzero.

Next, we look at a different class of matrices associated to a bipartite graph $G = (V^-,V^+,E)$, where the vertex sets $V^-$ and $V^+$ are given by
\begin{equation}
\label{vertices}
\begin{aligned}
V^- &= \{r_1,r_2,...,r_s\} \\
V^+ &= \{q_1,q_2,...,q_t\}.
\end{aligned}
\end{equation}
The pattern class $\mathcal{P}(G)$ of the bipartite graph $G$, with vertex sets $V^-$ and $V^+$ given by \eqref{vertices}, is defined as
\begin{equation}
\mathcal{P}(G) = \{ M \in \mathbb{R}^{t \times s} \: | \: M_{ij} \neq 0 \iff (r_j,q_i) \in E \}.
\end{equation}
Note that the cardinalities of $V^-$ and $V^+$ can differ, hence the matrices in the pattern class $\mathcal{P}(G)$ are not necessarily square.

\subsection{Subclass of distance-information preserving matrices}
In this subsection we investigate properties of the powers of matrices belonging to the qualitative class $\mathcal{Q}(G)$. The relevance of these properties will become apparent later on, when we provide a graph-theoretic condition for targeted controllability of systems defined on graphs. 

We first provide the following lemma, which states that if the distance between two nodes is greater than $k$, the corresponding element in $X^k$ is zero. 
\begin{lemma}
\label{Lemma1}
Consider a directed graph $G = (V,E)$, two distinct vertices $i,j \in V$, a matrix $X \in \mathcal{Q}(G)$ and a positive integer $k$. If $d(j,i) > k$, then $(X^k)_{ij} = 0$.
\end{lemma}
\begin{myproof}[Proof]
The proof follows easily by induction on $k$, and is therefore omitted. 
\end{myproof}

\noindent
Subsequently, we consider the class of matrices for which $(X^k)_{ij}$ is nonzero if the distance $d(j,i)$ is exactly equal to $k$. Such matrices are called distance-information preserving, more precisely:
\begin{dfn}
Consider a directed graph $G = (V,E)$. A matrix $X \in \mathcal{Q}(G)$ is called \textit{distance-information preserving} if for any two distinct vertices $i,j \in V$ we have that $d(j,i) = k$ implies $(X^k)_{ij} \neq 0$. 
\end{dfn}

Although the distance-information preserving property does not hold for all matrices $X \in \mathcal{Q}(G)$, it does hold for the adjacency and Laplacian matrices \cite{faults}.
Because these matrices are often used to describe network dynamics, distance-information preserving matrices form an important subclass of $\mathcal{Q}(G)$, which from now on will be denoted by $\mathcal{Q}_d(G)$. 

\subsection{Zero forcing sets}
In this section we review the notion of zero forcing. The reason for this is the correspondence between zero forcing sets and the sets of leaders rendering a system defined on a graph controllable. More on this will follow in the next subsection. 

For now, let $G = (V,E)$ be a directed graph with vertices colored either black or white. The color-change rule is defined in the following way: If $u \in V$ is a black vertex and exactly one out-neighbour $v \in V$ of $u$ is white, then change the color of $v$ to black \cite{hogben}.

When the color-change rule is applied to $u$ to change the color of $v$, we say $u$ forces $v$, and write $u \to v$.

Given a coloring of $G$, that is: given a set $C \subseteq V$ containing black vertices only, and a set $V \setminus C$ consisting of only white vertices, the derived set $D(C)$ is the set of black vertices obtained by applying the color-change rule until no more changes are possible \cite{hogben}. 

A zero forcing set for $G$ is a subset of vertices $Z \subseteq V$ such that if initially the vertices in $Z$ are colored black and the remaining vertices are colored white, then $D(Z) = V$. 

Finally, for a given zero forcing set, we can construct the derived set, listing the forces in the order in which they were performed. This list is called a chronological list of forces. Note that such a list does not have to be unique. 

\begin{example}
Consider the directed graph $G = (V,E)$ depicted in Figure \ref{fig:G1zfs} and let $C = \{2\}$.

\begin{minipage}{0.15\textwidth}
\begin{figure}[H]
\centering
\scalebox{0.8}{
\begin{tikzpicture}[scale=1]
		\node[draw, style=wball, label={90:$1$}] (1) at (0,0) {};
		\node[draw, style=bball, label={-90:$2$}] (2) at (0,-0.8) {};
		\node[draw, style=wball, label={90:$3$}] (3) at (0.8,0) {};
		\node[draw, style=wball, label={-90:$4$}] (4) at (0.8,-0.8) {};
		\node[draw, style=wball, label={-90:$5$}] (5) at (1.6,-0.8) {};
		\draw[-latex] (1) -- (3);
		\draw[-latex] (2) -- (4);
		\draw[-latex] (4) -- (5);
		\draw[-latex] (1) -- (4);
\end{tikzpicture}
}
\caption{Graph $G$.}
\label{fig:G1zfs}
\end{figure}
\end{minipage} \hfill
\begin{minipage}{0.15\textwidth}
\begin{figure}[H]
\centering
\scalebox{0.8}{
	\begin{tikzpicture}[scale=1]
		\node[draw, style=wball, label={90:$1$}] (1) at (0,0) {};
		\node[draw, style=bball, label={-90:$2$}] (2) at (0,-0.8) {};
		\node[draw, style=wball, label={90:$3$}] (3) at (0.8,0) {};
		\node[draw, style=bball, label={-90:$4$}] (4) at (0.8,-0.8) {};
		\node[draw, style=wball, label={-90:$5$}] (5) at (1.6,-0.8) {};
		\draw[-latex] (1) -- (3);
		\draw[-latex] (2) -- (4);
		\draw[-latex] (4) -- (5);
		\draw[-latex] (1) -- (4);
\end{tikzpicture}
}
\caption{Force $2 \to 4$.}
\label{fig:G2zfs}
\end{figure}
\end{minipage} \hfill
\begin{minipage}{0.15\textwidth}
\begin{figure}[H]
\centering
\scalebox{0.8}{
\begin{tikzpicture}[scale=1]
		\node[draw, style=wball, label={90:$1$}] (1) at (0,0) {};
		\node[draw, style=bball, label={-90:$2$}] (2) at (0,-0.8) {};
		\node[draw, style=wball, label={90:$3$}] (3) at (0.8,0) {};
		\node[draw, style=bball, label={-90:$4$}] (4) at (0.8,-0.8) {};
		\node[draw, style=bball, label={-90:$5$}] (5) at (1.6,-0.8) {};
		\draw[-latex] (1) -- (3);
		\draw[-latex] (2) -- (4);
		\draw[-latex] (4) -- (5);
		\draw[-latex] (1) -- (4);
\end{tikzpicture}
}
\caption{Force $4 \to 5$.}
\label{fig:G3zfs}
\end{figure}
\end{minipage}

\vspace{\belowdisplayskip}
\noindent
Note that vertex $2$ can force $4$, and subsequently node $4$ can force $5$. No further color changes can be made, so $D(C) = \{2,4,5\}$. As $D(C) \neq V$, $C$ is not a zero forcing set. However, suppose we choose $C = \{1,2\}$. In this case it is easy to see that we can color all vertices black, hence $C = \{1,2\}$ is a zero forcing set.
\end{example}

\subsection{Systems defined on graphs}
Consider a directed graph $G = (V,E)$, where the vertex set is given by $V = \{1,2,...,n\}$. Furthermore, let $V' = \{ v_1,v_2,...,v_r \} \subseteq V$ be a subset. The $n \times r$ matrix $P(V;V')$ is defined by 
\begin{equation}
P_{ij} = \left\{ \begin{matrix*}[l] 1 & \mbox{if } i = v_j \\ 0 & \mbox{otherwise. } \end{matrix*}\right.
\end{equation}

We now introduce the subset $V_L \subseteq V$ consisting of so-called \textit{leader nodes}, i.e. agents of the network to which an external control input is applied. The remaining nodes $V \setminus V_L$ are appropriately called followers. We consider finite-dimensional linear time-invariant systems of the form
\begin{equation}
\label{system}
\begin{aligned}
\dot{x}(t) = Xx(t) + Uu(t),
\end{aligned}
\end{equation}
where $x \in \mathbb{R}^n$ is the state and $u \in \mathbb{R}^m$ is the input of the system. Here $X \in \mathcal{Q}(G)$ and $U = P(V;V_L)$, for some leader set $V_L \subseteq V$. An important notion regarding systems of the form \eqref{system} is the notion of strong structural controllability. 
\begin{dfn}
\cite{zeroforce} A system of the form \eqref{system} is called \textit{strongly structurally controllable} if the pair $(X,U)$ is controllable for all $X \in \mathcal{Q}(G)$.
\end{dfn}

In the case that \eqref{system} is strongly structurally controllable we say $(G;V_L)$ is controllable, with a slight abuse of terminology. There is a one-to-one correspondence between strong structural controllability and zero forcing sets, as stated in the following theorem. 
\begin{theorem}
\label{Nimatheorem}
\cite{zeroforce} Let $G = (V,E)$ be a directed graph and let $V_L \subseteq V$ be a leader set. Then $(G;V_L)$ is controllable if and only if $V_L$ is a zero forcing set. 
\end{theorem}

In this paper, we are primarily interested in cases for which $(G;V_L)$ is not controllable. In such cases, we wonder whether we can control the state of a subset $V_T \subseteq V$ of nodes, called \textit{target nodes}. To this extent we specify an output equation $y(t) = Hx(t)$, which together with \eqref{system} yields the system
\begin{equation}
\label{fullsystem}
\begin{aligned}
\dot{x}(t) &= Xx(t) + Uu(t) \\
y(t) &= Hx(t),
\end{aligned}
\end{equation}
where $y \in \mathbb{R}^p$ is the output of the system consisting of the states of the target nodes, and $H = P^T(V;V_T)$. Note that the ability to control the states of all target nodes in $V_T$ is equivalent with the output controllability of system \eqref{fullsystem} \cite{Nima}. As the output of system \eqref{fullsystem} specifically consists of the states of the target nodes, we say \eqref{fullsystem} is targeted controllable if it is output controllable. 

Furthermore, system \eqref{fullsystem} is called strongly targeted controllable if $(X,U,H)$ is targeted controllable for all $X \in \mathcal{Q}(G)$ \cite{Nima}. In case \eqref{fullsystem} is strongly targeted controllable, we say $(G;V_L;V_T)$ is targeted controllable with respect to $\mathcal{Q}(G)$. The term `with respect to $\mathcal{Q}(G)$' clarifies the class of state matrices under consideration. This paper mainly considers strong targeted controllability with respect to $\mathcal{Q}_d(G)$. We conclude this section with well-known conditions for strong targeted controllability. Let $U = P(V;V_L)$ and $H = P^T(V;V_T)$ be the input and output matrices respectively, and define the reachable subspace 
$\big \langle X \: | \: \im U \big \rangle = \im U + X \im U + \cdots + X^{n-1} \im U$.
\begin{proposition}
\label{prop1}
\cite{Nima} The following statements are equivalent: \\[0.2cm]
1) $(G;V_L;V_T)$ is targeted controllable with respect to $\mathcal{Q}(G)$ \\
2) $\rank \begin{pmatrix} HU & HXU & \cdots & HX^{n-1}U \end{pmatrix} = p$ $\:\forall \: X \in \mathcal{Q}(G)$ \\
3) $H \big \langle X \: | \: \im U \big \rangle = \mathbb{R}^p$ $\:\forall \: X \in \mathcal{Q}(G)$ \\
4) $\ker H + \big \langle X \: | \: \im U \big \rangle = \mathbb{R}^n$ $\:\forall \: X  \in \mathcal{Q}(G)$.
\end{proposition}

\section{Problem statement}
\label{Problem Statement}
Strong targeted controllability with respect to $\mathcal{Q}(G)$ was studied in \cite{Nima}, and a sufficient graph-theoretic condition was provided. Motivated by the fact that $\mathcal{Q}_d(G)$ contains important network-related matrices like the adjacency and Laplacian matrices, we are interested in extending the results of \cite{Nima} to the class of distance-information preserving matrices $Q_d(G)$. More explicitly, the problem that we will investigate in this paper is given as follows.
\begin{problem}
\label{Problem1}
Given a directed graph $G = (V,E)$, a leader set $V_L \subseteq V$ and target set $V_T \subseteq V$, provide necessary and sufficient graph-theoretic conditions under which $(G;V_L;V_T)$ is targeted controllable with respect to $\mathcal{Q}_d(G)$.
\end{problem}
Such graph-theoretic conditions have a considerable advantage in their computational robustness compared to rank conditions, and aid in finding leader selection procedures. In addition, we are interested in a method to compute leader sets achieving targeted controllability. More precisely:
\begin{problem}
\label{Problem2}
Given a directed graph $G = (V,E)$ and target set $V_T \subseteq V$, compute a leader set $V_L \subseteq V$ of minimum cardinality such that $(G;V_L;V_T)$ is targeted controllable with respect to $\mathcal{Q}_d(G)$.
\end{problem}

\section{Main results}
\label{Main results}
Our main results are presented in this section. Firstly, in Section \ref{Sufficient condition} we  provide a sufficient graph-theoretic condition for strong targeted controllability with respect to $Q_d(G)$. Subsequently, in Section \ref{Sufficiently rich} we review the notion of sufficient richness of subclasses, and prove that the subclass $\mathcal{Q}_d(G)$ is sufficiently rich. This result allows us to establish a necessary condition for strong targeted controllability, which is presented in Section \ref{Necessary condition}. Finally, in Section \ref{algorithmsection} we show Problem \ref{Problem2} is NP-hard and provide a heuristic leader selection algorithm to determine leader sets achieving targeted controllability. 

\subsection{Sufficient condition for targeted controllability}
\label{Sufficient condition}
This section discusses a sufficient graph-theoretic condition for strong targeted controllability. We first introduce some notions that will become useful later on. 

Consider a directed graph $G = (V,E)$ with a leader set $V_L$ and target set $V_T$. The derived set of $V_L$ is given by $D(V_L)$. Furthermore, let $V_S \subseteq V \setminus D(V_L)$ be a subset. We partition the set $V_S$ according to the distance of its nodes with respect to $D(V_L)$, that is
\begin{equation}
\label{partition}
V_S = V_1 \cup V_2 \cup \dots \cup V_d,
\end{equation}
where for $j \in V_S$ we have $j \in V_i$ if and only if $d(D(V_L),j) = i$ for $i = 1,2,...,d$. Moreover, we define $\check{V}_i$ and $\hat{V}_i$ to be the sets of vertices in $V_S$ of distance respectively less than $i$ and greater than $i$ with respect to $D(V_L)$. More precisely:
\begin{equation}
\begin{aligned}
\check{V}_i &:= V_1 \cup ... \cup V_{i-1} \text{ for } i = 2,...,d \\
\hat{V}_i &:= V_{i+1} \cup ... \cup V_{d} \text{ for } i = 1,...,d-1
\end{aligned}
\end{equation} 
By convention $\check{V}_1 = \varnothing$ and $\hat{V}_d = \varnothing$.
With each of the sets $V_1, V_2,...,V_d$ we associate a bipartite graph $G_i = (D(V_L),V_i, E_i)$, where for $j \in D(V_L)$ and $k \in V_i$ we have $(j,k) \in E_i$ if and only if $d(j,k) = i$ in the network graph $G$. 

\begin{example}
\label{bipartitegraphexample}
We consider the network graph $G = (V,E)$ as depicted in Figure \ref{fig:network}. The set of leaders is $V_L = \{1,2\}$, which implies that $D(V_L) = \{1,2,3\}$. 

\begin{minipage}{0.23\textwidth}
			\begin{figure}[H]
				\centering
				\scalebox{0.8}{
					\begin{tikzpicture}[scale=1]
					\tikzset{dot/.style={circle,fill=#1,inner sep=0,minimum size=4pt}}
					\node[draw, style=bball, label={90:$1$}] (1) at (0,0) {};
					\node[draw, style=bball, label={-90:$2$}] (2) at (0,-1.5) {};
					\node[draw, style=wball, label={90:$3$}] (3) at (1.5,0) {};
					\node[draw, style=wball, label={-90:$4$}] (4) at (1.5,-1.5) {};
					\node[draw, style=wball, label={-90:$5$}] (5) at (1.5,-2.5) {};
					\node[draw, style=wball, label={90:$9$}] (6) at (3,0) {};
					\node[draw, style=wball, label={90:$7$}] (7) at (4.5,0) {};
					\node[draw, style=wball, label={-90:$6$}] (8) at (3,-1.5) {};
					\node[draw, style=wball, label={90:$10$}] (9) at (4.5,1) {};
					\node[draw, style=wball, label={-90:$8$}] (10) at (4.5,-1.5) {};
					\draw[-latex] (1) -- (3);
					\draw[-latex] (2) -- (4);
					\draw[-latex] (2) -- (5);
					\draw[-latex] (3) -- (6);
					\draw[-latex] (6) -- (7);
					\draw[-latex] (4) -- (8);
					\draw[-latex] (4) -- (6);
					\draw[-latex] (8) -- (10);
					\draw[-latex] (6) -- (9);
					\draw[-latex] (3) to [bend right=30] (4);
					\draw[-latex] (4) to [bend right=30] (3);
					\draw[-latex] (6) to [bend right=30] (8);
					\draw[-latex] (8) to [bend right=30] (6);
					\end{tikzpicture}
				}
				\caption{Graph $G$ with $V_L = \{1,2\}$.}
				\label{fig:network}
			\end{figure}
		\end{minipage} \hfill
		\begin{minipage}{0.23\textwidth}
			\begin{figure}[H]
				\centering
				\scalebox{0.8}{
					\begin{tikzpicture}[scale=1]
					\tikzset{dot/.style={circle,fill=#1,inner sep=0,minimum size=4pt}}
					\node[draw, style=bball, label={90:$1$}] (1) at (0,0) {};
					\node[draw, style=bball, label={-90:$2$}] (2) at (0,-1.5) {};
					\node[draw, style=bball, label={90:$3$}] (3) at (1.5,0) {};
					\node[draw, style=wball, label={-90:$4$}] (4) at (1.5,-1.5) {};
					\node[draw, style=wball, label={-90:$5$}] (5) at (1.5,-2.5) {};
					\node[draw, style=wball, label={90:$9$}] (6) at (3,0) {};
					\node[draw, style=wball, label={90:$7$}] (7) at (4.5,0) {};
					\node[draw, style=wball, label={-90:$6$}] (8) at (3,-1.5) {};
					\node[draw, style=wball, label={90:$10$}] (9) at (4.5,1) {};
					\node[draw, style=wball, label={-90:$8$}] (10) at (4.5,-1.5) {};
					\draw[-latex] (1) -- (3);
					\draw[-latex] (2) -- (4);
					\draw[-latex] (2) -- (5);
					\draw[-latex] (3) -- (6);
					\draw[-latex] (6) -- (7);
					\draw[-latex] (4) -- (8);
					\draw[-latex] (4) -- (6);
					\draw[-latex] (8) -- (10);
					\draw[-latex] (6) -- (9);
					\draw[-latex] (3) to [bend right=30] (4);
					\draw[-latex] (4) to [bend right=30] (3);
					\draw[-latex] (6) to [bend right=30] (8);
					\draw[-latex] (8) to [bend right=30] (6);
					\end{tikzpicture}
				}
				\caption{$D(V_L) = \{1,2,3\}$.}
				\label{fig:Exchange}
			\end{figure}
		\end{minipage}

\vspace{\belowdisplayskip}

In this example, we define the subset $V_S \subseteq V \setminus D(V_L)$ as $V_S := \{4,5,6,7,8\}$. Note that $V_S$ can be partitioned according to the distance of its nodes with respect to $D(V_L)$ as $V_S = V_1 \cup V_2 \cup V_3$, where $V_1 = \{4,5\}$, $V_2 = \{6,7\}$ and $V_3 = \{8\}$. The bipartite graphs $G_1$, $G_2$ and $G_3$ are given in Figures \ref{fig:G1},\ref{fig:G2} and \ref{fig:G3} respectively.

\begin{minipage}{0.15\textwidth}
\begin{figure}[H]
\centering
\scalebox{0.8}{
\begin{tikzpicture}[scale=1]
		\tikzset{dot/.style={circle,fill=#1,inner sep=0,minimum size=4pt}}
		\node[draw, style=wball, minimum size=0.4cm, label={0:$4$}] (4) at (2,0) {};
		\node[draw, style=wball, minimum size=0.4cm, label={0:$5$}] (5) at (2,-1) {};
		\node[draw, style=wball, minimum size=0.4cm, label={180:$1$}] (1) at (0,0.5) {};
		\node[draw, style=wball, minimum size=0.4cm, label={180:$2$}] (2) at (0,-0.5) {};
		\node[draw, style=wball, minimum size=0.4cm, label={180:$3$}] (3) at (0,-1.5) {};
		\draw[-latex] (2) -- (4);
		\draw[-latex] (2) -- (5);		
		\draw[-latex] (3) -- (4);
		\end{tikzpicture}
}
\caption{Graph $G_1$.}
\label{fig:G1}
\end{figure}
\end{minipage} \hfill
\begin{minipage}{0.15\textwidth}
\begin{figure}[H]
\centering
\scalebox{0.8}{
	\begin{tikzpicture}[scale=1]
		\node[draw, style=wball, minimum size=0.4cm, label={0:$6$}] (6) at (2,0) {};
		\node[draw, style=wball, minimum size=0.4cm, label={0:$7$}] (7) at (2,-1) {};
		\node[draw, style=wball, minimum size=0.4cm, label={180:$1$}] (1) at (0,0.5) {};
		\node[draw, style=wball, minimum size=0.4cm, label={180:$2$}] (2) at (0,-0.5) {};
		\node[draw, style=wball, minimum size=0.4cm, label={180:$3$}] (3) at (0,-1.5) {};
		\draw[-latex] (2) -- (6);
		\draw[-latex] (3) -- (6);		
		\draw[-latex] (3) -- (7);
	\end{tikzpicture}
}
\caption{Graph $G_2$.}
\label{fig:G2}
\end{figure}
\end{minipage} \hfill
\begin{minipage}{0.15\textwidth}
\begin{figure}[H]
\centering
\scalebox{0.8}{
\begin{tikzpicture}[scale=1]
		\node[draw, style=wball, minimum size=0.4cm, label={0:$8$}] (8) at (2,-0.5) {};
		\node[draw, style=wball, minimum size=0.4cm, label={180:$1$}] (1) at (0,0.5) {};
		\node[draw, style=wball, minimum size=0.4cm, label={180:$2$}] (2) at (0,-0.5) {};
		\node[draw, style=wball, minimum size=0.4cm, label={180:$3$}] (3) at (0,-1.5) {};
		\draw[-latex] (2) -- (8);
		\draw[-latex] (3) -- (8);		
		\end{tikzpicture}
}
\caption{Graph $G_3$.}
\label{fig:G3}
\end{figure}
\end{minipage}
\end{example}
\vspace{\belowdisplayskip}
\noindent
The main result presented in this section is given in Theorem \ref{main}. This statement provides a sufficient graph-theoretic condition for targeted controllability of $(G;V_L;V_T)$ with respect to $\mathcal{Q}_d(G)$.
\begin{theorem}
\label{main}
Consider a directed graph $G = (V,E)$, with leader set $V_L \subseteq V$ and target set $V_T \subseteq V$. Let $V_T \setminus D(V_L)$ be partitioned as in \eqref{partition}, and assume $D(V_L)$ is a zero forcing set in $G_i = (D(V_L),V_i,E_i)$ for $i = 1,2,...,d$. Then $(G;V_L;V_T)$ is targeted controllable with respect to $\mathcal{Q}_d(G)$.
\end{theorem}
The `k-walk theory' for target control, described in \cite{barabasitarget} is just a special case of Theorem \ref{main}. Indeed, in the single-leader case, the condition of Theorem \ref{main} reduces to the condition that no pair of target nodes has the same distance with respect to the leader. However, it is worth mentioning that Theorem \ref{main} holds for general directed graphs and multiple leaders, while the results of \cite{barabasitarget} are only applicable to directed tree networks in the case that the leader set is singleton. Furthermore, note that Theorem \ref{main} significantly improves the known condition for strong targeted controllability given in \cite{Nima} for the class $\mathcal{Q}_d(G)$. In Theorem \ref{main} target nodes with arbitrary distance with respect to the derived set are allowed, while the main result Theorem VI.6 of \cite{Nima} is restricted to target nodes of distance one with respect to $D(V_L)$. Before proving Theorem \ref{main}, we provide an illustrative example and two auxiliary lemmas.
\begin{example}
\label{Examplesufficientcondition}
Once again, consider the network graph depicted in Figure \ref{fig:network}, with leader set $V_L = \{1,2\}$ and assume the target set is given by $V_T = \{ 1,2,...,8 \}$. The goal of this example is to prove that $(G;V_L;V_T)$ is targeted controllable with respect to $\mathcal{Q}_d(G)$. 

Note that $V_S := V_T \setminus D(V_L)$ is given by $V_S = \{ 4,5,6,7,8 \}$, which is partitioned according to \eqref{partition} as $V_S = V_1 \cup V_2 \cup V_3$, where $V_1 = \{4,5\}$, $V_2 = \{6,7\}$ and $V_3 = \{8\}$. The graphs $G_1$, $G_2$ and $G_3$ have been computed in Example \ref{bipartitegraphexample}. Note that $D(V_L) = \{1,2,3\}$ is a zero forcing set in all three graphs (see Figures \ref{fig:G1}, \ref{fig:G2} and \ref{fig:G3}). We conclude by Theorem \ref{main} that $(G;V_L;V_T)$ is targeted controllable with respect to $\mathcal{Q}_d(G)$.
\end{example}

\begin{lemma}
\label{derivedsetlemma}
Consider a directed graph $G = (V,E)$ with leader set $V_L \subseteq V$ and target set $V_T \subseteq V$. Let $\mathcal{Q}_s(G) \subseteq \mathcal{Q}(G)$ be any subclass. Then $(G;V_L;V_T)$ is targeted controllable with respect to $Q_s(G)$ if and only if $(G;D(V_L);V_T)$ is targeted controllable with respect to $Q_s(G)$.
\end{lemma}
\begin{myproof}
Let $U = P(V;V_L)$ index the leader set $V_L$ and $W = P(V;D(V_L))$ index the derived set of $V_L$. We have that $(G;V_L;V_T)$ is targeted controllable with respect to $Q_s(G)$ if and only if
\begin{equation}
\label{condition1}
H \big \langle X \: | \: \im U \big \rangle = \mathbb{R}^p \text{ for all } X \in Q_s(G),
\end{equation}
However, as $\big \langle X \: | \: \im U \big \rangle = \big \langle X \: | \: \im W \big \rangle$ for any $X \in \mathcal{Q}(G)$ (see Lemma VI.2 of \cite{Nima}), \eqref{condition1} holds if and only if
\begin{equation}
H \big \langle X \: | \: \im W \big \rangle = \mathbb{R}^p \text{ for all } X \in Q_s(G).
\end{equation}
We conclude that $(G;V_L;V_T)$ is targeted controllable with respect to $Q_s(G)$ if and only if $(G; D(V_L); V_T)$ is targeted controllable with respect to $Q_s(G)$.
\end{myproof}

\begin{lemma}
\label{ZFSbipartitelemma}
Let $G = (V^-,V^+,E)$ be a bipartite graph and assume $V^-$ is a zero forcing set in $G$. Then all matrices $M \in \mathcal{P}(G)$ have full row rank. 
\end{lemma}
\begin{myproof}
Note that forces of the form $u \to v$, where $u,v \in V^+$ are not possible, as $G$ is a bipartite graph. Relabel the nodes of $V^-$ and $V^+$ such that a chronological list of forces is given by $u_i \to v_i$, where $u_i \in V^-$ and $v_i \in V^+$ for $i = 1,2,...,|V^+|$. Let $M \in \mathcal{P}(G)$ be a matrix in the pattern class of $G$. Note that the element $M_{ii}$ is nonzero, as $u_i \to v_i$. Furthermore, $M_{ji}$ is zero for all $j > i$. The latter follows from the fact that $u_i$ would not be able to force $v_i$ if there was an arc $(u_i,v_j) \in E$. We conclude that the columns $1,2,...,|V^+|$ of $M$ are linearly independent, hence $M$ has full row rank. 
\end{myproof}

\begin{myproof}[Proof of Theorem \ref{main}]
Let $D(V_L) = \{ 1,2,...,m \}$, and assume without loss of generality that the matrix $U$ has the form (see Lemma \ref{derivedsetlemma}):
\begin{equation}
U = \begin{pmatrix}
I_{m\times m} & 0_{m \times (n-m)}
\end{pmatrix}^T .
\end{equation}
Furthermore, we let $V_S := V_T \setminus D(V_L)$ be given by $\{ m+1,m+2,...,p \}$, where the vertices are ordered in non-decreasing distance with respect to $D(V_L)$. Partition $V_S$ according to the distance of its nodes with respect to $D(V_L)$ as
\begin{equation}
V_S = V_1 \cup V_2 \cup \dots \cup V_d,
\end{equation}
where  for $j \in V_S$ we have $j \in V_i$ if and only if $d(D(V_L),j) = i$ for $i = 1,2,...,d$. Finally, assume the target set $V_T$ contains all nodes in the derived set $D(V_L)$. This implies that the matrix $H$ is of the form
\begin{equation}
H = \begin{pmatrix}
I_{p\times p} & 0_{p \times (n-p)} 
\end{pmatrix}.
\end{equation}
Note that by the structure of $H$ and $U$, the matrix $HX^iU$ is simply the $p \times m$ upper left corner submatrix of $X^i$. We now claim that $HX^iU$ can be written as follows.
\begin{equation}
\label{HXU}
HX^iU = \begin{pmatrix}
\Lambda_i \\[0.2cm]
  M_i  \\[0.2cm]
 0_i
\end{pmatrix},
\end{equation}
where $M_i \in \mathcal{P}(G_i)$ is a $|V_i| \times m$ matrix in the pattern class of $G_i$, $\Lambda_i$ is an $\left(m + |\check{V}_i|\right) \times m$ matrix containing elements of lesser interest, and $0_i$ is a zero matrix of dimension $|\hat{V}_i| \times m$.

We proceed as follows: first we prove that the bottom submatrix of \eqref{HXU} contains zeros only, secondly we prove that $M_i \in \mathcal{P}(G_i)$. From this, we conclude that equation \eqref{HXU} holds.

Note that for $k \in D(V_L)$ and $j \in \hat{V}_i$, we have $d(k,j) > i$ and by Lemma \ref{Lemma1} it follows that $(X^i)_{jk} = 0$. As $D(V_L) = \{1,2,...,m\}$, this means that the bottom $|\hat{V}_i| \times m$ submatrix of $HX^iU$ is a zero matrix.

Subsequently, we want to prove that $M_i$, the middle block of \eqref{HXU}, is an element of the pattern class $\in \mathcal{P}(G_i)$. Note that the $j$th row of $M_i$ corresponds to the element $l := m + |\check{V}_i| + j \in V_i$.

Suppose $(M_i)_{jk} \neq 0$ for a $k \in \{1,2,...,m\}$ and $j \in \{1,2,...,|V_i|\}$.
As $M_i$ is a submatrix of $HX^iU$, this implies $(HX^iU)_{lk} \neq 0$. Recall that $HX^iU$ is the $p \times m$ upper left corner submatrix of $X^i$, therefore it holds that $(X^i)_{lk} \neq 0$.  Note that for the vertices $k \in D(V_L)$ and $l \in V_i$ we have $d(k,l) \geq i$ by the partition of $V_S$. However, as $(X^i)_{lk} \neq 0$ it follows from Lemma \ref{Lemma1} that $d(k,l) = i$. Therefore, by the definition of $G_i$, there is an arc $(k,l) \in E_i$. 

Conversely, suppose there is an arc $(k,l) \in E_i$ for $l \in V_i$ and $k \in D(V_L)$. This implies $d(k,l) = i$ in the network graph $G$. By the distance-information preserving property of $X$ we consequently have $(X^i)_{lk} \neq 0$. We conclude that $(M_i)_{jk} \neq 0$ and hence $M_i \in \mathcal{P}(G_i)$. This implies that equation \eqref{HXU} holds, 
We compute the first $dm$ columns of the output controllability matrix $\begin{pmatrix} HU & HXU & HX^2U & \dots & HX^dU \end{pmatrix}$ as follows:
\begin{equation}
\label{HUHXU}
\begin{pmatrix}
I & \ast & \ast & \dots & \ast & \ast \\
0 & M_1 & \ast & \dots & \ast & \ast \\
0 & 0 & M_2 & \ddots & \vdots & \vdots \\
0 & 0 & 0 & \ddots & \ast & \ast \\
\vdots & \vdots & \vdots & \ddots & M_{d-1} & \ast  \\
0 & 0 & 0 & \dots & 0 & M_d
\end{pmatrix},
\end{equation}
where zeros denote zero matrices and asterisks denote matrices of less interest. As $D(V_L)$ is a zero forcing set in $G_i$ for $i = 1,2,...,d$, the matrices $M_1,M_2,...,M_d$ have full row rank by Lemma \ref{ZFSbipartitelemma}. We conclude that the matrix \eqref{HUHXU} has full row rank, and consequently $(G;V_L;V_T)$ is targeted controllable with respect to $\mathcal{Q}_d(G)$. 
\end{myproof}

Note that the condition given in Theorem \ref{main} is sufficient, but not necessary. One can verify that the graph $G = (V,E)$ with leader set $V_L = \{1\}$ and target set $V_T = \{2,3\}$ depicted in Figure \ref{fig:examplenotnec} is an example of a graph for which $(G;V_L;V_T)$ is targeted controllable with respect to $\mathcal{Q}_d(G)$. However, this graph does not satisfy the conditions of Theorem \ref{main}.

\begin{figure}[h]
\begin{minipage}[b]{0.23\textwidth}
\centering
\scalebox{0.8}{
\begin{tikzpicture}[scale=1]
		\tikzset{dot/.style={circle,fill=#1,inner sep=0,minimum size=4pt}}
		\node[draw, style=bball, label={90:$1$}] (1) at (0,0) {};
		\node[draw, style=wball, label={90:$2$}] (2) at (2,1) {};
		\node[draw, style=wball, label={90:$3$}] (3) at (2,-1) {};
		\node[draw, style=wball, label={90:$4$}] (4) at (4,0) {};
		\draw[-latex] (1) -- (2);
		\draw[-latex] (1) -- (3);
		\draw[-latex] (2) -- (4);
		\draw[-latex] (4) -- (3);
		\end{tikzpicture}
}
\end{minipage} \hfill
\begin{minipage}[b]{0.23\textwidth}
\centering
\scalebox{0.8}{
\begin{tikzpicture}[scale=1]
		\tikzset{dot/.style={circle,fill=#1,inner sep=0,minimum size=4pt}}
		\node[draw, style=bball, label={90:$1$}] (1) at (0,0) {};
		\node[draw, style=wball, label={90:$2$}] (3) at (2,1) {};
		\node[draw, style=wball, label={90:$4$}] (5) at (4,1) {};
		\node[draw, style=wball, label={90:$3$}] (2) at (2,-1) {};
		\node[draw, style=wball, label={90:$5$}] (4) at (4,-1) {};
		\draw[-latex] (1) -- (2);
		\draw[-latex] (2) -- (4);
		\draw[-latex] (1) -- (3);
		\draw[-latex] (3) -- (5);
		\end{tikzpicture}
}
\end{minipage} \hfill

 \begin{minipage}[t]{.5\linewidth}
 \centering
        \caption{Theorem \ref{main} not necessary.}
		\label{fig:examplenotnec}
    \end{minipage}%
    \hfill%
    \begin{minipage}[t]{.5\linewidth}
    \centering
        \caption{Theorem \ref{necessaryconditiontc} not sufficient.}
		\label{fig:notsuf2}
    \end{minipage}%
\end{figure}

\subsection{Sufficient richness of $Q_d(G)$}
\label{Sufficiently rich}
The notion of sufficient richness of a qualitative subclass was introduced in \cite{zeroforce}. We provide an equivalent definition as follows.
\begin{dfn}
Let $G = (V,E)$ be a directed graph with leader set $V_L \subseteq V$. A subclass $Q_s(G) \subseteq \mathcal{Q}(G)$ is called \textit{sufficiently rich} if $(G;V_L)$ is controllable with respect to $Q_s(G)$ implies $(G;V_L)$ is controllable with respect to $Q(G)$. 
\end{dfn}

The following geometric characterization of sufficient richness is proven in \cite{zeroforce}.
\begin{proposition}
\label{sufrichprop}
A qualitative subclass $Q_s(G) \subseteq \mathcal{Q}(G)$ is sufficiently rich if for all $z \in \mathbb{R}^n$ and $X \in \mathcal{Q}(G)$ satisfying $z^T X = 0$, there exists an $X'\in Q_s(G)$ such that $z^T X' = 0$.
\end{proposition}

The goal of this section is to prove that the qualitative subclass of distance-information preserving matrices is sufficiently rich. This result will be used later on, when we provide a necessary condition for targeted controllability with respect to $\mathcal{Q}_d(G)$. First however, we state two auxiliary lemmas which will be the building blocks to prove the sufficient richness of $\mathcal{Q}_d(G)$. 
\begin{lemma}
\label{propnonzero}
Consider $q$ nonzero multivariate polynomials $p_i(x)$, where $i = 1,2,...,q$ and $x \in \mathbb{R}^n$. There exists an $\bar{x} \in \mathbb{R}^n$ such that $p_i(\bar{x}) \neq 0$ for $i = 1,2,...,q$.   
\end{lemma}
\begin{myproof}
The proof follows immediately from continuity of polynomials and is omitted.
\end{myproof}
\noindent
\textbf{Remark:} Without loss of generality, we can assume that the point $\bar{x} \in \mathbb{R}^n$ has only nonzero coordinates. Indeed, if $p_i(\bar{x}) \neq 0$ for $i = 1,2,...,q$, there exists an open ball $B(\bar{x})$ around $\bar{x}$ in which $p_i(x) \neq 0$ for $i = 1,2,...,q$. Obviously, this open ball contains a point with the aforementioned property.

\begin{lemma}
\label{proppol}
Let $X \in \mathcal{Q}(G)$ and $D = \text{diag }(d_1,d_2,...,d_n)$ be a matrix with variable diagonal entries. If $d(i,j) = k$ for distinct vertices $i$ and $j$, then $((XD)^k)_{ji}$ is a nonzero polynomial in the variables $d_1,d_2,...,d_n$. 
\end{lemma}

\begin{myproof}
Note that $((XD)^k)_{ji}$ is given by 
\begin{equation*}
\begin{aligned}
\sum\limits_{i_1 = 1}^n \sum\limits_{i_2 = 1}^n \dotsb \sum\limits_{i_{k-1} = 1}^n (XD)_{i_1,i} (XD)_{i_2,i_1} \dotsb (XD)_{j,i_{k-1}},
\end{aligned}
\end{equation*}
which equals
\begin{equation}
\label{sumoverpaths}
\begin{aligned}
\sum\limits_{i_1 = 1}^n \sum\limits_{i_2 = 1}^n \dotsb \sum\limits_{i_{k-1} = 1}^n  d_{i} X_{i_1,i} \cdot d_{i_1} X_{i_2,i_1} \dotsb d_{i_{k-1}} X_{j,i_{k-1}}.
\end{aligned}
\end{equation}
Since the distance $d(i,j)$ is equal to $k$, there exists at least one path of length $k$ from $i$ to $j$, which we denote by $(i,i_1),(i_1,i_2),...,(i_{k-1},j)$. It follows that the corresponding elements of the matrix $X$, i.e. the elements $X_{i_1,i}, X_{i_2,i_1}, X_{j,i_{k-1}}$ are nonzero. Therefore, the term 
\begin{equation}
\label{nonvanishingterm}
d_{i} X_{i_1,i} \cdot d_{i_1} X_{i_2,i_1} \dotsb d_{i_{k-1}} X_{j,i_{k-1}}
\end{equation}
is nonzero (as a function of $d_i,d_{i_1},d_{i_2},...,d_{i_{k-1}}$). Furthermore, this combination of $k$ diagonal elements is unique in the sense that there does not exist another summand on the right-hand side of \eqref{sumoverpaths} with exactly the same elements. This implies that the term \eqref{nonvanishingterm} does not vanish (as a polynomial). We conclude that $((XD)^k)_{ji}$ is a nonzero polynomial function in the variables $d_1,d_2,...,d_n$.
\end{myproof}

\begin{theorem}
The subclass $\mathcal{Q}_d(G)$ is sufficiently rich. 
\end{theorem}
\begin{myproof}
Given a matrix $X \in \mathcal{Q}(G)$, using Lemmas \ref{propnonzero} and \ref{proppol}, we first prove there exists a diagonal matrix $\bar{D}$ with nonzero diagonal components such that $X\bar{D} \in \mathcal{Q}_d(G)$. From this we will conclude $\mathcal{Q}_d(G)$ is sufficiently rich. 

Let $D = \text{diag }(d_1,d_2,...,d_n)$ be a matrix with variable diagonal entries. We define $p_{ij} := ((XD)^{d(i,j)})_{ji}$ for distinct $i,j = 1,2,...,n$. By Lemma \ref{proppol} we have that $p_{ij}(d_1,d_2,...,d_n)$ is a nonzero polynomial in the variables $d_1,d_2,...,d_n$. Moreover, Lemma \ref{propnonzero} states the existence of nonzero real constants $\bar{d_1},\bar{d_2},...,\bar{d_n}$ such that
\begin{equation}
p_{ij}(\bar{d_1},\bar{d_2},...,\bar{d_n}) \neq 0 \text{ for distinct } i,j = 1,2...,n.
\end{equation}
Therefore, the choice $\bar{D} = \text{diag }(\bar{d_1},\bar{d_2},...,\bar{d_n})$ implies $X\bar{D} \in \mathcal{Q}_d(G)$. Let $z \in \mathbb{R}^n$ be a vector such that $z^T X = 0$ for an $X \in \mathcal{Q}(G)$. The choice of $X' = X\bar{D}$ yields a matrix $X' \in \mathcal{Q}_d(G)$ for which $z^T X' = 0$. By Proposition \ref{sufrichprop} it follows that $\mathcal{Q}_d(G)$ is sufficiently rich. 
\end{myproof}
\subsection{Necessary condition for targeted controllability}
\label{Necessary condition}
In addition to the previously established sufficient condition for targeted controllability, we give a necessary graph-theoretic condition for targeted controllability in Theorem \ref{necessaryconditiontc}. 
\begin{theorem}
\label{necessaryconditiontc}
Let $G = (V,E)$ be a directed graph with leader set $V_L \subseteq V$ and target set $V_T \subseteq V$. If $(G;V_L;V_T)$ is targeted controllable with respect to $\mathcal{Q}_d(G)$ then $V_L \cup (V \setminus V_T)$ is a zero forcing set in $G$.
\end{theorem}
\begin{myproof}
Assume without loss of generality that $V_L \cap V_T = \varnothing$. Hence, $V_L \cup (V \setminus V_T) = V \setminus V_T$. Partition the vertex set into $V_L$, $V \setminus (V_L \cup V_T)$ and $V_T$. Accordingly, the input and output matrices $U = P(V;V_L)$ and $H = P^T(V;V_T)$ satisfy
\begin{equation}
U = \begin{pmatrix}
I & 0 & 0
\end{pmatrix}^T \qquad 
H = \begin{pmatrix}
0 & 0 & I
\end{pmatrix}
\end{equation}
Note that $\ker H = \im R$, where $R = P(V;(V \setminus V_T))$ is given by
\begin{equation}
R = \begin{pmatrix}
I & 0 & 0 \\
0 & I & 0
\end{pmatrix}^T.
\end{equation}
Since for all $X \in \mathcal{Q}_d(G)$ we have
\begin{equation}
\begin{aligned}
\ker H + \big \langle X \: | \: \im U \big \rangle &= \mathbb{R}^n,
\end{aligned}
\end{equation}
equivalently,
\begin{equation}
\begin{aligned}
\im R + \big \langle X \: | \: \im U \big \rangle &= \mathbb{R}^n,
\end{aligned}
\end{equation}
we obtain
\begin{equation}
\label{kern}
\begin{aligned}
\big \langle X \: | \: \im \begin{pmatrix} \: U & R \: \end{pmatrix} \big \rangle &= \mathbb{R}^n.
\end{aligned}
\end{equation}
As $\im U \subseteq \im R$, \eqref{kern} implies $\big \langle X \: | \: \im R \: \big \rangle = \mathbb{R}^n$ for all $X \in \mathcal{Q}_d(G)$, equivalently, the pair $(X,R)$ is controllable for all $X \in \mathcal{Q}_d(G)$. However, by sufficient richness of  $\mathcal{Q}_d(G)$, it follows that $(X,R)$ is controllable for all $X \in \mathcal{Q}(G)$.
We conclude from Theorem \ref{Nimatheorem} that $V \setminus V_T$ is a zero forcing set. 
\end{myproof}
\noindent
\begin{example}
Consider the directed graph $G = (V,E)$ with leader set $V_L = \{1,2\}$ and target set $V_T = \{1,2,...,8\}$ as depicted in Figure \ref{fig:network}. We know from Example \ref{Examplesufficientcondition} that $(G;V_L;V_T)$ is targeted controllable with respect to $\mathcal{Q}_d(G)$. The set $V_L \cup (V \setminus V_T) = \{1,2,9,10\}$ is colored black in Figure \ref{fig:examplegraphnec}. Indeed, $V_L \cup (V \setminus V_T)$ is a zero forcing set in $G$. A possible chronological list of forces is: $1 \to 3$, $3 \to 4$, $2 \to 5$, $4 \to 6$, $6 \to 8$ and $9 \to 7$. 
\begin{figure}[h]
\begin{minipage}[b]{0.48\textwidth}
\centering
\scalebox{0.8}{
\begin{tikzpicture}[scale=1]
		\tikzset{dot/.style={circle,fill=#1,inner sep=0,minimum size=4pt}}
		\node[draw, style=bball, label={90:$1$}] (1) at (0,0) {};
		\node[draw, style=bball, label={-90:$2$}] (2) at (0,-1.5) {};
		\node[draw, style=wball, label={90:$3$}] (3) at (1.5,0) {};
		\node[draw, style=wball, label={-90:$4$}] (4) at (1.5,-1.5) {};
		\node[draw, style=wball, label={-90:$5$}] (5) at (1.5,-2.5) {};
		\node[draw, style=bball, label={90:$9$}] (6) at (3,0) {};
		\node[draw, style=wball, label={90:$7$}] (7) at (4.5,0) {};
		\node[draw, style=wball, label={-90:$6$}] (8) at (3,-1.5) {};
		\node[draw, style=bball, label={90:$10$}] (9) at (4.5,1) {};
		\node[draw, style=wball, label={-90:$8$}] (10) at (4.5,-1.5) {};
		\draw[-latex] (1) -- (3);
		\draw[-latex] (2) -- (4);
		\draw[-latex] (2) -- (5);
		\draw[-latex] (3) -- (6);
		\draw[-latex] (6) -- (7);
		\draw[-latex] (4) -- (8);
		\draw[-latex] (4) -- (6);
		\draw[-latex] (8) -- (10);
		\draw[-latex] (6) -- (9);
		\draw[-latex] (3) to [bend right=30] (4);
		\draw[-latex] (4) to [bend right=30] (3);
		\draw[-latex] (6) to [bend right=30] (8);
		\draw[-latex] (8) to [bend right=30] (6);
		\node[draw, style=wball, label={90:$1$}] (11) at (6,0) {};
		\node[draw, style=wball, label={-90:$2$}] (12) at (6,-1.5) {};
		\node[draw, style=wball, label={90:$3$}] (13) at (8,0) {};
		\node[draw, style=wball, label={-90:$4$}] (14) at (8,-1.5) {};
		\node[draw, style=wball, label={90:$5$}] (15) at (10,1.5) {};
		\node[draw, style=wball, label={90:$6$}] (16) at (10,0) {};
		\node[draw, style=wball, label={-90:$7$}] (17) at (10,-1.5) {};
		\draw[-latex] (11) -- (13);
		\draw[-latex] (12) -- (13);
		\draw[-latex] (12) -- (14);
		\draw[-latex] (13) -- (15);
		\draw[-latex] (13) -- (16);
		\draw[-latex] (13) -- (17);
		\draw[-latex] (14) -- (17);
		\end{tikzpicture}
}
\end{minipage} \hfill
\begin{minipage}[b]{0\textwidth}
\centering
\scalebox{0.8}{
\begin{tikzpicture}[scale=1]
		\node[] (1) at (0,0) {};
		\end{tikzpicture}
}
\end{minipage} \hfill

 \begin{minipage}[t]{.5\linewidth}
 \centering
        \caption{Zero forcing set.}
        \label{fig:examplegraphnec}
    \end{minipage}%
    \hfill%
    \begin{minipage}[t]{.5\linewidth}
    \centering
        \caption{Example of Algorithm 1.}
        \label{fig:examplefase1}
    \end{minipage}%
\end{figure}
\end{example}
The condition provided in Theorem \ref{necessaryconditiontc} is necessary for targeted controllability, but not sufficient. To prove this fact, consider the directed graph with leader set $V_L = \{1\}$ and target set $V_T = \{4,5\}$ given in Figure \ref{fig:notsuf2}. It can be shown that $(G;V_L;V_T)$ is not targeted controllable with respect to $\mathcal{Q}_d(G)$, even though $V_L \cup (V \setminus V_T) = \{1,2,3\}$ is a zero forcing set.

So far, we have provided a necessary and a sufficient topological condition for targeted controllability. However, given a network graph with target set, it is not clear how to choose leaders achieving target control. Hence, in the following section we focus on a leader selection algorithm. 

\subsection{Leader selection algorithm}
\label{algorithmsection}
The problem addressed in this section is as follows: given a directed graph $G = (V,E)$ with target set $V_T \subseteq V$, find a leader set $V_L \subseteq V$ of minimum cardinality such that $(G;V_L;V_T)$ is targeted controllable with respect to $\mathcal{Q}_d(G)$. Such a leader set is called a minimum leader set. 

As the problem of finding a minimum zero forcing set is NP-hard \cite{hardness}, the problem of finding a minimum leader set $V_L$ that achieves controllability of $(G;V_L)$ is NP-hard. Since this controllability problem can be seen as a special case of target control, where the entire vertex set is regarded as target set, we therefore conclude that determining a minimum leader set such that $(G;V_L;V_T)$ is targeted controllable with respect to $\mathcal{Q}_d(G)$ is NP-hard. 

It is for this reason we propose a heuristic approach to compute a (minimum) leader set that achieves targeted controllability. The algorithm consists of two phases. Firstly, we identify a set of nodes in the graph $G$ from which all target nodes can be reached. These nodes are taken as leaders. Secondly, this set of leaders is extended to achieve targeted controllability. 

To explain the first phase of the algorithm, we introduce some notation. First of all, we define the notion of \textit{root set}. 
\begin{dfn}
Consider a directed graph $G = (V,E)$ and a target set $V_T \subseteq V$. A subset $V_R \subseteq V$ is called a root set of $V_T$ if for any $v \in V_T$ there exists a vertex $u \in V_R$ such that $d(u,v) < \infty$.
\end{dfn}

A root set of $V_T$ of minimum cardinality is called a minimum root set of $V_T$. Note that the cardinality of a minimum root set of $V_T$ is a lower bound on the minimum number of leaders rendering $(G;V_L;V_T)$ targeted controllable. Indeed, it is easy to see that if there are no paths from any of the leader nodes to a target node, the graph is not targeted controllable. The first step of the proposed algorithm is to compute the minimum root set of $V_T$. Let the vertex and target sets be given by $V = \{ 1,2,...,n \}$ and $V_T = \{ v_1,v_2,...,v_p \} \subseteq V$ respectively. Furthermore, define a matrix $A \in \mathbb{R}^{p \times n}$ in the following way. For $j \in V$ and $v_i \in V_T$ let 
\begin{equation}
\label{matrixA}
A_{ij} := \left\{
                \begin{array}{ll}
                1 & \text{if } d(j,v_i) < \infty \\
                0 & \text{otherwise}
                \end{array}
              \right.
\end{equation}
That is: the matrix $A$ contains zeros and ones only, where coefficients with value one indicate the existence of a path between the corresponding vertices. Finding a minimum root set of $V_T$ boils down to finding a binary vector $x \in \mathbb{R}^n$ with minimum number of ones such that $Ax \geq \mathbbm{1}$, where the inequality is defined element-wise and $\mathbbm{1}$ denotes the vector of all ones. In this vector $x$, coefficients with value one correspond to elements in the root set of $V_T$. It is for this reason we can formulate the minimum root set problem as a binary integer linear program
\begin{equation}
\label{intprog}
\begin{aligned}
	\text{minimize } &\mathbbm{1}^T x \\
	\text{subject to } &Ax \geq \mathbbm{1} \\
	\text{and } &x \in \{0,1\}^n.
\end{aligned}
\end{equation}
Linear programs of this form can be solved using software like CPLEX or Matlab. For very large-scale problems one might resort to heuristic methods (see e.g. \cite{BDD}). In the following example we illustrate how the minimum root set problem can be regarded as a binary integer linear program.
\begin{example}
\label{examplefase1}
Consider the directed graph $G = (V,E)$ with target set $V_T = \{1,4,5,6,7\}$ depicted in Figure \ref{fig:examplefase1}. The goal of this example is to find a minimum root set for $V_T$. The matrix $A$, as defined in \eqref{matrixA}, is given by
\begin{equation}
A = \begin{pmatrix}
1 & 0 & 0 & 0 & 0 & 0 & 0 \\
0 & 1 & 0 & 1 & 0 & 0 & 0 \\
1 & 1 & 1 & 0 & 1 & 0 & 0 \\
1 & 1 & 1 & 0 & 0 & 1 & 0 \\
1 & 1 & 1 & 1 & 0 & 0 & 1
\end{pmatrix}.
\end{equation}
Note that $x = \begin{pmatrix}
1 & 0 & 0 & 1 & 0 & 0 & 0
\end{pmatrix}^T$ satisfies the inequality $Ax \geq \mathbbm{1}$ and the constraint $x \in \{0,1\}^7$. Furthermore, the vector $x$ minimizes $\mathbbm{1}^T x$ under these constraints. This can be seen in the following way: there is no column of $A$ in which all elements equal 1, hence there is no vector $x$ with a single one such that $Ax \geq \mathbbm{1}$ is satisfied. Therefore, $x$ solves the binary integer linear program \eqref{intprog}, from which we conclude that the choice $V_R = \{1,4\}$ yields a minimum root set for $V_T$. Indeed, observe in Figure \ref{fig:examplefase1} that we can reach all nodes in the target set starting from the nodes $1$ and $4$. It is worth mentioning that the choice of minimum root set is not unique: the set $\{1,2\}$ is also a minimum root set for $V_T$.
\end{example}

In general, the minimum root set $V_R$ of $V_T$ does not guarantee targeted controllability of $(G;V_R;V_T)$ with respect to $\mathcal{Q}_d(G)$. For instance, it can be shown for the graph $G$ and target set $V_T$ of Example \ref{examplefase1} that the leader set $V_L = \{1,4\}$ does not render $(G;V_L;V_T)$ targeted controllable with respect to $\mathcal{Q}_d(G)$. Hence, we propose a greedy approach to extend the minimum root set of $V_T$ to a leader set that does achieve targeted controllability. 

Recall from Theorem \ref{main} that $(G;V_L;V_T)$ is targeted controllable with respect to $\mathcal{Q}_d(G)$ if $D(V_L)$ is a zero forcing set in the bipartite graphs $G_i = (D(V_L),V_i,E_i)$ for $i = 1,2,...,d$, where $V_i \subseteq V_T$ is the set of target nodes having distance $i$ from $D(V_L)$. Given an initial set of leaders $V_L$, we compute its derived set $D(V_L)$ and  verify whether we can force all nodes in the bipartite graphs $G_i$ for $i = 1,2,...,d$. Suppose that in the bipartite graph $G_k$ the set $V_k$ cannot be forced by $D(V_L)$ for a $k \in \{1,2,...,d\}$. Let $V_k = \{v_{1},v_{2},...,v_{l}\}$, and suppose $v_{i}$ is the first vertex in $V_k$ that cannot be forced. Then we choose $v_i$ as a leader. Consequently, we have extended our leader set $V_L$ to $V_L \cup \{v_i\}$. With the extended leader set we can repeat the procedure, until the leaders render the graph targeted controllable. This idea is captured more formally in the following leader selection algorithm. One should recognize the two phases of leader selection: firstly, a minimum root set is computed. Subsequently, the minimum root set is greedily extended to a leader set achieving targeted controllability.

\noindent
\hrulefill \\
 \textbf{Algorithm 1:} Leader Selection Procedure 
 
\noindent 
\hrulefill

{\setlength{\abovedisplayskip}{0pt}%
\begin{flalign*}
\begin{array}{ll}
\textbf{Input: } & \text{Directed graph } G = (V,E); \\
 & \text{Target set } V_T \subseteq V; \\
\textbf{Output: } & \text{Leader set } V_L \subseteq V \text{ achieving target control};
\end{array}
\end{flalign*}%
}%

\hspace{0.4cm} Let $V_L = \varnothing$; 

\hspace{0.4cm} Compute matrix $A$, given in \eqref{matrixA}; 

\hspace{0.4cm} Find a solution $x$ for the linear program \eqref{intprog}; 

\hspace{0.4cm} \textbf{for} $i = 1$ to $n$;

\hspace{0.8cm} \textbf{if} $x_i = 1$;

\hspace{1.2cm} $V_L$ = $V_L \cup i$;

\hspace{0.8cm} \textbf{end}

\hspace{0.4cm} \textbf{end}

\hspace{0.4cm} Compute $D(V_L)$; 

\hspace{0.4cm} Set $i = 1$; 

\hspace{0.4cm} \textbf{repeat} 

\hspace{0.8cm} Compute $V_i$ and $G_i = (D(V_L),V_i,E_i)$;

\hspace{0.8cm} \textbf{if} $D(V_L)$ forces $V_i$ in $G_i$;

\hspace{1.2cm} $i = i+1$;

\hspace{0.8cm} \textbf{else} 

\hspace{1.2cm} Let $v$ be the first unforced vertex in $V_i$;

\hspace{1.2cm} Set $V_L = V_L \cup v$;

\hspace{1.2cm} Compute $D(V_L)$;

\hspace{1.2cm} $i = 1$;

\hspace{0.8cm} \textbf{end} 

\hspace{0.4cm} \textbf{until} $d(D(V_L),v) < i$ for all $v \in V_T$; 

\textbf{return } $V_L$.

\noindent \hrulefill

\begin{figure*}
\centering
\scalebox{1}{
\begin{tikzpicture}[scale=0.95]
		
\node[draw, style=wball, label={90:$1$}] (1) at (0,1) {};
\node[draw, circle, minimum size=1.3cm, color = red, thick] at (1.6,0.8) {};
\node[draw, style=wball, label={90:$2$}] (2) at (1.6,0.8) {};
\node[draw, circle, minimum size=1.3cm, color = red, thick] at (2,-1.2) {};
\node[draw, style=wball, label={90:$3$}] (3) at (2,-1.2){};
\node[draw, style=wball, label={90:$4$}] (4) at (3.4,0.15) {};
\node[draw, style=wball, label={90:$5$}] (5) at (5.4,0.1) {};
\node[draw, circle, minimum size=1.3cm, color = red, thick] at (7.4,0) {};
\node[draw, style=wball, label={90:$6$}] (6) at (7.4,0) {};
\node[draw, style=wball, label={90:$7$}] (7) at (10,0.6) {};
\node[draw, circle, minimum size=1.3cm, color = red, thick] at (9.6,-0.5) {};
\node[draw, style=wball, label={-90:$8$}] (8) at (9.6,-0.5) {};

\node[draw, style=wball, label={90:$9$}] (9) at (12.3,1.6) {};
\node[draw, circle, minimum size=1.3cm, color = red, thick] at (11.5,-0.1) {};
\node[draw, style=wball, label={90:$10$}] (10) at (11.5,-0.1) {};
\node[draw, style=wball, label={-90:$11$}] (11) at (11.4,-1.2) {};
\node[draw, style=wball, label={-90:$12$}] (12) at (12.8,-2) {};

\node[draw, circle, minimum size=1.3cm, color = red, thick] at (13.8,2.7) {};
\node[draw, style=wball, label={90:$13$}] (13) at (13.8,2.7) {};
\node[draw, style=wball, label={90:$14$}] (14) at (14.7,1.4) {};
\node[draw, circle, minimum size=1.3cm, color = red, thick] at (13.5,0) {};
\node[draw, style=wball, label={90:$15$}] (15) at (13.5,0) {};
\node[draw, circle, minimum size=1.3cm, color = red, thick] at (14.5,-2.3) {};
\node[draw, style=wball, label={-90:$16$}] (16) at (14.5,-2.3) {};

\node[draw, circle, minimum size=1.3cm, color = red, thick] at (16.2,2.3) {};
\node[draw, style=wball, label={90:$17$}] (17) at (16.2,2.3) {};
\node[draw, style=wball, label={-90:$18$}] (18) at (15.5,-0.7) {};
\node[draw, style=wball, label={90:$19$}] (19) at (17.5,1.4) {};
\node[draw, circle, minimum size=1.3cm, color = red, thick] at (17,0.3) {};
\node[draw, style=wball, label={-90:$20$}] (20) at (17,0.3) {};

\draw[-latex] (1) -- (2);
\draw[-latex] (4) -- (2);
\draw[-latex] (4) -- (3);
\draw[-latex] (4) -- (5);
\draw[-latex] (5) -- (6);
\draw[-latex] (7) -- (6);
\draw[-latex] (8) -- (6);
\draw[-latex] (7) -- (8);	
\draw[-latex] (9) -- (7);	
\draw[-latex] (7) -- (10);	
\draw[-latex] (10) -- (9);	
\draw[-latex] (8) -- (11);	
\draw[-latex] (10) -- (12);	
\draw[-latex] (13) -- (9);	
\draw[-latex] (13) -- (14);	
\draw[-latex] (9) -- (14);	
\draw[-latex] (15) -- (14);	
\draw[-latex] (10) -- (15);	
\draw[-latex] (15) -- (11);	
\draw[-latex] (15) -- (16);	
\draw[-latex] (16) -- (12);	
\draw[-latex] (13) -- (17);	
\draw[-latex] (14) -- (17);	
\draw[-latex] (14) -- (19);	
\draw[-latex] (17) -- (20);	
\draw[-latex] (19) -- (20);	
\draw[-latex] (18) -- (20);	
\draw[-latex] (18) -- (16);	
\draw[-latex] (15) -- (18);	
\draw[-latex] (11) -- (18);	
	
\end{tikzpicture}
}
\caption{Directed graph $G = (V,E)$ with encircled target nodes $V_T = \{2,3,6,8,10,13,15,16,17,20\}$.}
\label{fig:fullscaleexample}
\end{figure*}
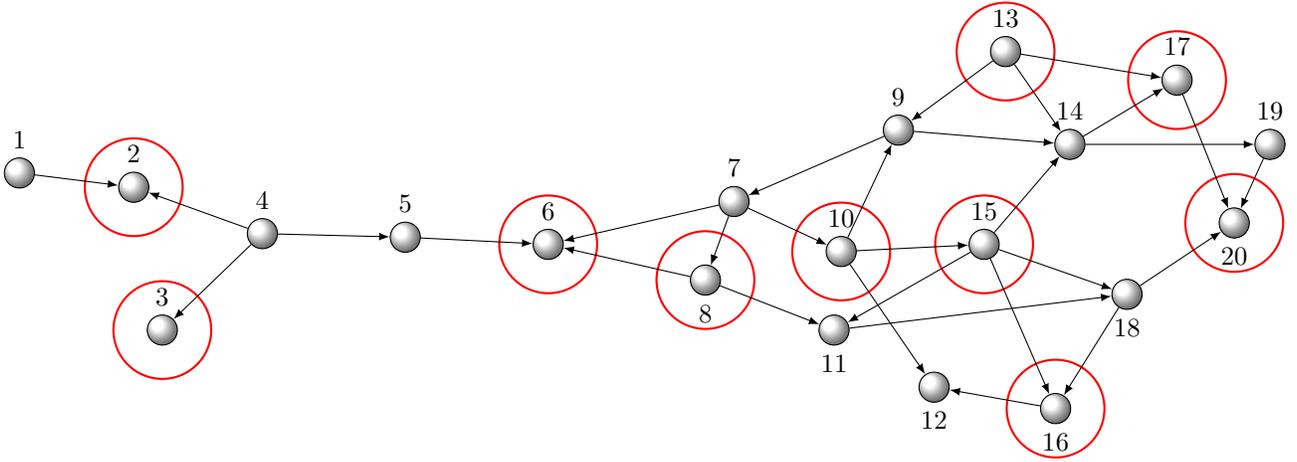
\begin{example}
 Consider the directed graph given in Figure \ref{fig:fullscaleexample}, with target set $V_T = \{2,3,6,8,10,13,15,16,17,20\}$. The goal of this example is to compute a leader set $V_L$ such that $(G;V_L;V_T)$ is targeted controllable with respect to $\mathcal{Q}_d(G)$. 
The first step of Algorithm 1 is to compute the matrix $A$, defined in \eqref{matrixA}. For this example, $A$ is given as follows.
\begin{equation*}
\resizebox{.95\hsize}{!}{$
A = \begin{pmatrix}
1 & 1 & 0 & 1 & 0 & 0 & 0 & 0 & 0 & 0 & 0 & 0 & 0 & 0 & 0 & 0 & 0 & 0 & 0 & 0 \\
0 & 0 & 1 & 1 & 0 & 0 & 0 & 0 & 0 & 0 & 0 & 0 & 0 & 0 & 0 & 0 & 0 & 0 & 0 & 0 \\
0 & 0 & 0 & 1 & 1 & 1 & 1 & 1 & 1 & 1 & 0 & 0 & 1 & 0 & 0 & 0 & 0 & 0 & 0 & 0 \\
0 & 0 & 0 & 0 & 0 & 0 & 1 & 1 & 1 & 1 & 0 & 0 & 1 & 0 & 0 & 0 & 0 & 0 & 0 & 0 \\
0 & 0 & 0 & 0 & 0 & 0 & 1 & 0 & 1 & 1 & 0 & 0 & 1 & 0 & 0 & 0 & 0 & 0 & 0 & 0 \\
0 & 0 & 0 & 0 & 0 & 0 & 0 & 0 & 0 & 0 & 0 & 0 & 1 & 0 & 0 & 0 & 0 & 0 & 0 & 0 \\
0 & 0 & 0 & 0 & 0 & 0 & 1 & 0 & 1 & 1 & 0 & 0 & 1 & 0 & 1 & 0 & 0 & 0 & 0 & 0 \\
0 & 0 & 0 & 0 & 0 & 0 & 1 & 1 & 1 & 1 & 1 & 0 & 1 & 0 & 1 & 1 & 0 & 1 & 0 & 0 \\
0 & 0 & 0 & 0 & 0 & 0 & 1 & 0 & 1 & 1 & 0 & 0 & 1 & 1 & 1 & 0 & 1 & 0 & 0 & 0 \\
0 & 0 & 0 & 0 & 0 & 0 & 1 & 1 & 1 & 1 & 1 & 0 & 1 & 1 & 1 & 0 & 1 & 1 & 1 & 1
\end{pmatrix}$}
\end{equation*}
Using the Matlab function {\tt intlinprog}, we find the optimal solution 
\begin{equation*}
\resizebox{.95\hsize}{!}{$
x = \begin{pmatrix}
0 & 0 & 0 & 1 & 0 & 0 & 0 & 0 & 0 & 0 & 0 & 0 & 1 & 0 & 0 & 0 & 0 & 0 & 0 & 0
\end{pmatrix}^T$}
\end{equation*}
to the binary linear program \eqref{intprog}. Hence, a minimum root set for $V_T$ is given by $\{4,13\}$. Following Algorithm 1, we define our initial leader set $V_L = \{4,13\}$. As nodes $4$ and $13$ both have three out-neighbours, the derived set of $V_L$ is simply given by $D(V_L) = \{4,13\}$. The next step of the algorithm is to compute the first bipartite graph $G_1 = (D(V_L),V_1, E_1)$, which we display in Figure \ref{fig:firstbip}.

\begin{minipage}{0.22\textwidth}
\begin{figure}[H]
\centering
\scalebox{0.8}{
\begin{tikzpicture}[scale=1]
		\node[draw, style=wball, minimum size=0.4cm, label={0:$2$}] (2) at (2,1) {};
		\node[draw, style=wball, minimum size=0.4cm, label={0:$3$}] (3) at (2,0) {};
		\node[draw, style=wball, minimum size=0.4cm, label={0:$17$}] (17) at (2,-1) {};
		\node[draw, style=bball, minimum size=0.4cm, label={180:$4$}] (4) at (0,0.5) {};
		\node[draw, style=bball, minimum size=0.4cm, label={180:$13$}] (13) at (0,-0.5){};
		\draw[-latex] (13) -- (17);
		\draw[-latex] (4) -- (2);
		\draw[-latex] (4) -- (3);
\end{tikzpicture}
}
\caption{Graph $G_1$ for derived set $D(V_L) = \{4,13\}$.}
\label{fig:firstbip}
\end{figure}
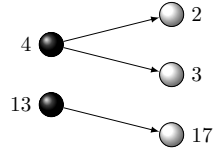
\end{minipage} \hfill
\begin{minipage}{0.22\textwidth}
\begin{figure}[H]
\centering
\scalebox{0.8}{
\begin{tikzpicture}[scale=1]
		\node[draw, style=wball, minimum size=0.4cm, label={0:$3$}] (3) at (2,0.5) {};
		\node[draw, style=wball, minimum size=0.4cm, label={0:$17$}] (17) at (2,-.5) {};
		\node[draw, style=bball, minimum size=0.4cm, label={180:$2$}] (2) at (0,1) {};
		\node[draw, style=bball, minimum size=0.4cm, label={180:$4$}] (4) at (0,0) {};
		\node[draw, style=bball, minimum size=0.4cm, label={180:$13$}] (13) at (0,-1) {};
		\draw[-latex] (13) -- (17);
		\draw[-latex] (4) -- (3);
\end{tikzpicture}
}
\caption{Graph $G_1$ for derived set $D(V_L) = \{2,4,13\}$.}
\label{fig:firstbip2}
\end{figure}
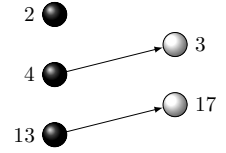
\end{minipage}

\vspace{\belowdisplayskip}
\noindent
Observe that the nodes $2$ and $3$ cannot be forced, hence we choose node $2$ as additional leader. The process now repeats itself, we redefine $V_L = \{2,4,13\}$ and compute $D(V_L) = \{2,4,13\}$. Furthermore, for this leader set, the graph $G_1 = (D(V_L),V_1,E_1)$ is given in Figure \ref{fig:firstbip2}. In this case, the set $V_1 = \{3,17\}$ of nodes having distance one with respect to $D(V_L)$ is forced. Therefore, we continue with the second bipartite graph $G_2 = (D(V_L),V_2,E_2)$, given in Figure \ref{fig:secondbip2}. 

\begin{minipage}{0.22\textwidth}
\begin{figure}[H]
\centering
\scalebox{0.8}{
\begin{tikzpicture}[scale=1]
		\node[draw, style=wball, minimum size=0.4cm, label={0:$6$}] (3) at (2,0.5) {};
		\node[draw, style=wball, minimum size=0.4cm, label={0:$20$}] (17) at (2,-.5) {};
		\node[draw, style=bball, minimum size=0.4cm, label={180:$2$}] (2) at (0,1) {};
		\node[draw, style=bball, minimum size=0.4cm, label={180:$4$}] (4) at (0,0) {};
		\node[draw, style=bball, minimum size=0.4cm, label={180:$13$}] (13) at (0,-1) {};
		\draw[-latex] (13) -- (17);
		\draw[-latex] (4) -- (3);
\end{tikzpicture}
}
\caption{Graph $G_2$ for derived set $D(V_L) = \{2,4,13\}$.}
\label{fig:secondbip2}
\end{figure}
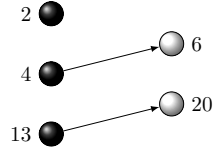
\end{minipage} \hfill
\begin{minipage}{0.22\textwidth}
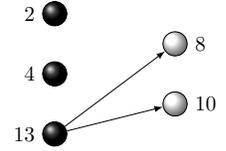
\begin{figure}[H]
\centering
\scalebox{0.8}{
\begin{tikzpicture}[scale=1]
		\node[draw, style=wball, minimum size=0.4cm, label={0:$8$}] (3) at (2,0.5) {};
		\node[draw, style=wball, minimum size=0.4cm, label={0:$10$}] (17) at (2,-.5) {};
		\node[draw, style=bball, minimum size=0.4cm, label={180:$2$}] (2) at (0,1) {};
		\node[draw, style=bball, minimum size=0.4cm, label={180:$4$}] (4) at (0,0) {};
		\node[draw, style=bball, minimum size=0.4cm, label={180:$13$}] (13) at (0,-1) {};
		\draw[-latex] (13) -- (17);
		\draw[-latex] (13) -- (3);
\end{tikzpicture}
}
\caption{Graph $G_3$ for derived set $D(V_L) = \{2,4,13\}$.}
\label{fig:thirdbip}
\end{figure}
\end{minipage}

\vspace{\belowdisplayskip}
\noindent
The set $V_2$ is forced by $D(V_L)$ in the graph $G_2$, hence we continue to investigate the third bipartite graph consisting of nodes having distance three with respect to $D(V_L)$. This graph is displayed in Figure \ref{fig:thirdbip}. As neither node $8$ nor $10$ can be forced, we add node $8$ to the leader set. In other words, we redefine $V_L = \{2,4,8,13\}$. Furthermore, the derived set of $V_L$ is given by $D(V_L) = \{2,4,8,13\}$. As we adapted the derived set, we have to recalculate the bipartite graphs $G_1$, $G_2$, $G_3$ and $G_4$ (see Figures \ref{fig:111}, \ref{fig:222}, \ref{fig:333} and \ref{fig:444}). 

\begin{minipage}{0.22\textwidth}
\begin{figure}[H]
\centering
\scalebox{0.8}{
\begin{tikzpicture}[scale=1]
		\node[draw, style=wball, minimum size=0.4cm, label={0:$3$}] (3) at (2,1) {};
		\node[draw, style=wball, minimum size=0.4cm, label={0:$6$}] (6) at (2,0) {};
		\node[draw, style=wball, minimum size=0.4cm, label={0:$17$}] (17) at (2,-1) {};
		\node[draw, style=bball, minimum size=0.4cm, label={180:$2$}] (2) at (0,1.5) {};
		\node[draw, style=bball, minimum size=0.4cm, label={180:$4$}] (4) at (0,0.5) {};
		\node[draw, style=bball, minimum size=0.4cm, label={180:$8$}] (8) at (0,-0.5) {};
		\node[draw, style=bball, minimum size=0.4cm, label={180:$13$}] (13) at (0,-1.5) {};
		\draw[-latex] (4) -- (3);
		\draw[-latex] (13) -- (17);
		\draw[-latex] (8) -- (6);
\end{tikzpicture}
}
\caption{Graph $G_1$ for derived set $D(V_L) = \{2,4,8,13\}$.}
\label{fig:111}
\end{figure}
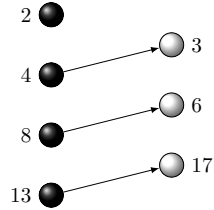
\end{minipage} \hfill
\begin{minipage}{0.22\textwidth}
\begin{figure}[H]
\centering
\scalebox{0.8}{
\begin{tikzpicture}[scale=1]
		\node[draw, style=wball, minimum size=0.4cm, label={0:$20$}] (20) at (2,0) {};
		\node[draw, style=bball, minimum size=0.4cm, label={180:$2$}] (2) at (0,1.5) {};
		\node[draw, style=bball, minimum size=0.4cm, label={180:$4$}] (4) at (0,0.5) {};
		\node[draw, style=bball, minimum size=0.4cm, label={180:$8$}] (8) at (0,-0.5) {};
		\node[draw, style=bball, minimum size=0.4cm, label={180:$13$}] (13) at (0,-1.5) {};
		\draw[-latex] (13) -- (20);
\end{tikzpicture}
}
\caption{Graph $G_2$ for derived set $D(V_L) = \{2,4,8,13\}$.}
\label{fig:222}
\end{figure}
\end{minipage}

\begin{minipage}{0.22\textwidth}
\begin{figure}[H]
\centering
\scalebox{0.8}{
\begin{tikzpicture}[scale=1]
		\node[draw, style=wball, minimum size=0.4cm, label={0:$10$}] (10) at (2,.5) {};
		\node[draw, style=wball, minimum size=0.4cm, label={0:$16$}] (16) at (2,-.5) {};
		\node[draw, style=bball, minimum size=0.4cm, label={180:$2$}] (2) at (0,1.5) {};
		\node[draw, style=bball, minimum size=0.4cm, label={180:$4$}] (4) at (0,0.5) {};
		\node[draw, style=bball, minimum size=0.4cm, label={180:$8$}] (8) at (0,-0.5) {};
		\node[draw, style=bball, minimum size=0.4cm, label={180:$13$}] (13) at (0,-1.5) {};
		\draw[-latex] (8) -- (16);
		\draw[-latex] (13) -- (10);
\end{tikzpicture}
}
\caption{Graph $G_3$ for derived set $D(V_L) = \{2,4,8,13\}$.}
\label{fig:333}
\end{figure}
\end{minipage} \hfill
\begin{minipage}{0.22\textwidth}
\begin{figure}[H]
\centering
\scalebox{0.8}{
\begin{tikzpicture}[scale=1]
		\node[draw, style=wball, minimum size=0.4cm, label={0:$15$}] (15) at (2,0) {};
		\node[draw, style=bball, minimum size=0.4cm, label={180:$2$}] (2) at (0,1.5) {};
		\node[draw, style=bball, minimum size=0.4cm, label={180:$4$}] (4) at (0,0.5) {};
		\node[draw, style=bball, minimum size=0.4cm, label={180:$8$}] (8) at (0,-0.5) {};
		\node[draw, style=bball, minimum size=0.4cm, label={180:$13$}] (13) at (0,-1.5) {};
		\draw[-latex] (13) -- (15);
\end{tikzpicture}
}
\caption{Graph $G_4$ for derived set $D(V_L) = \{2,4,8,13\}$.}
\label{fig:444}
\end{figure}
\end{minipage}

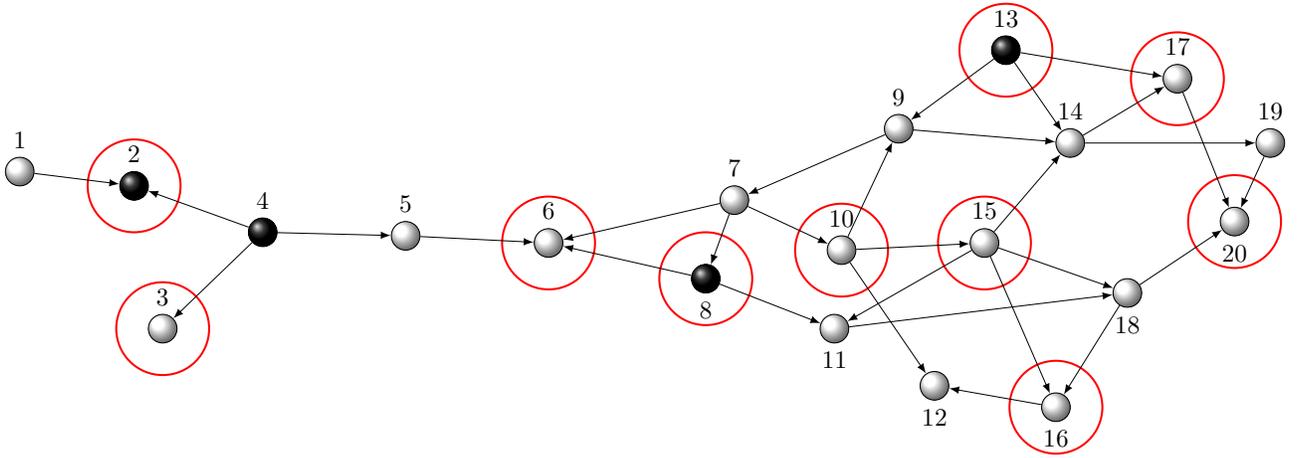
\begin{figure*}
\centering
\scalebox{0.95}{
\begin{tikzpicture}[scale=1]
		
\node[draw, style=wball, label={90:$1$}] (1) at (0,1) {};
\node[draw, circle, minimum size=1.3cm, color = red, thick] at (1.6,0.8) {};
\node[draw, style=bball, label={90:$2$}] (2) at (1.6,0.8) {};
\node[draw, circle, minimum size=1.3cm, color = red, thick] at (2,-1.2) {};
\node[draw, style=wball, label={90:$3$}] (3) at (2,-1.2){};
\node[draw, style=bball, label={90:$4$}] (4) at (3.4,0.15) {};
\node[draw, style=wball, label={90:$5$}] (5) at (5.4,0.1) {};
\node[draw, circle, minimum size=1.3cm, color = red, thick] at (7.4,0) {};
\node[draw, style=wball, label={90:$6$}] (6) at (7.4,0) {};
\node[draw, style=wball, label={90:$7$}] (7) at (10,0.6) {};
\node[draw, circle, minimum size=1.3cm, color = red, thick] at (9.6,-0.5) {};
\node[draw, style=bball, label={-90:$8$}] (8) at (9.6,-0.5) {};

\node[draw, style=wball, label={90:$9$}] (9) at (12.3,1.6) {};
\node[draw, circle, minimum size=1.3cm, color = red, thick] at (11.5,-0.1) {};
\node[draw, style=wball, label={90:$10$}] (10) at (11.5,-0.1) {};
\node[draw, style=wball, label={-90:$11$}] (11) at (11.4,-1.2) {};
\node[draw, style=wball, label={-90:$12$}] (12) at (12.8,-2) {};

\node[draw, circle, minimum size=1.3cm, color = red, thick] at (13.8,2.7) {};
\node[draw, style=bball, label={90:$13$}] (13) at (13.8,2.7) {};
\node[draw, style=wball, label={90:$14$}] (14) at (14.7,1.4) {};
\node[draw, circle, minimum size=1.3cm, color = red, thick] at (13.5,0) {};
\node[draw, style=wball, label={90:$15$}] (15) at (13.5,0) {};
\node[draw, circle, minimum size=1.3cm, color = red, thick] at (14.5,-2.3) {};
\node[draw, style=wball, label={-90:$16$}] (16) at (14.5,-2.3) {};

\node[draw, circle, minimum size=1.3cm, color = red, thick] at (16.2,2.3) {};
\node[draw, style=wball, label={90:$17$}] (17) at (16.2,2.3) {};
\node[draw, style=wball, label={-90:$18$}] (18) at (15.5,-0.7) {};
\node[draw, style=wball, label={90:$19$}] (19) at (17.5,1.4) {};
\node[draw, circle, minimum size=1.3cm, color = red, thick] at (17,0.3) {};
\node[draw, style=wball, label={-90:$20$}] (20) at (17,0.3) {};

\draw[-latex] (1) -- (2);
\draw[-latex] (4) -- (2);
\draw[-latex] (4) -- (3);
\draw[-latex] (4) -- (5);
\draw[-latex] (5) -- (6);
\draw[-latex] (7) -- (6);
\draw[-latex] (8) -- (6);
\draw[-latex] (7) -- (8);	
\draw[-latex] (9) -- (7);	
\draw[-latex] (7) -- (10);	
\draw[-latex] (10) -- (9);	
\draw[-latex] (8) -- (11);	
\draw[-latex] (10) -- (12);	
\draw[-latex] (13) -- (9);	
\draw[-latex] (13) -- (14);	
\draw[-latex] (9) -- (14);	
\draw[-latex] (15) -- (14);	
\draw[-latex] (10) -- (15);	
\draw[-latex] (15) -- (11);	
\draw[-latex] (15) -- (16);	
\draw[-latex] (16) -- (12);	
\draw[-latex] (13) -- (17);	
\draw[-latex] (14) -- (17);	
\draw[-latex] (14) -- (19);	
\draw[-latex] (17) -- (20);	
\draw[-latex] (19) -- (20);	
\draw[-latex] (18) -- (20);	
\draw[-latex] (18) -- (16);	
\draw[-latex] (15) -- (18);	
\draw[-latex] (11) -- (18);	
	
\end{tikzpicture}
}
\caption{Network graph $G = (V,E)$ with encircled target nodes $V_T = \{2,3,6,8,10,13,15,16,17,20\}$, and leader set $V_L = \{2,4,8,13\}$ colored black.}
\label{fig:fullscaleexample2}
\end{figure*}

\vspace{\belowdisplayskip}
\noindent
Note that in this case $D(V_L)$ is a zero forcing set in all four bipartite graphs. Furthermore, since $d(D(V_L),v) < 5$ for all $v \in V_T$, Algorithm 1 returns the leader set $V_L = \{2,4,8,13\}$. This choice of leader set guarantees that $(G;V_L;V_T)$ is targeted controllable with respect to $\mathcal{Q}_d(G)$. For the sake of clarity, we display the network graph in Figure \ref{fig:fullscaleexample2}, where the leader nodes are colored black, and the target nodes are encircled. 
\end{example}
It is worth mentioning that Algorithm 1 can also be applied to compute a leader set $V_L$ such that $(G;V_L)$ is controllable. Indeed, as full control can be regarded as a special case of target control where the entire vertex set $V$ is chosen as target set, Algorithm 1 can simply be applied to a directed graph $G = (V,E)$ with target set $V$. For example, we can retrieve the results on minimum leader sets found in \cite{zeroforce} for cycle and complete graphs, using Algorithm 1 and the fact that $\mathcal{Q}_d(G)$ is sufficiently rich. That is: for cycle and complete graphs, Algorithm 1 returns leader sets of respectively 2 and $n-1$ leaders, which is in agreement with \cite{zeroforce}.

\section{Conclusions}
\label{Conclusions}
In this paper, strong targeted controllability for the class of distance-information preserving matrices has been discussed. We have provided a sufficient graph-theoretic condition for strong targeted controllability, expressed in terms of zero-forcing sets of particular distance-related bipartite graphs. We have shown that this result significantly improves the known sufficient topological condition \cite{Nima} for strong targeted controllability of the class of distance-information preserving matrices. Furthermore, our result contains the `k-walk theory' \cite{barabasitarget} as a special case.

Motivated by the observation that the aforementioned sufficient condition is not a one-to-one correspondence, we provided a necessary topological condition for strong targeted controllability. This condition was proved using the fact that the subclass of distance-information preserving matrices is sufficiently rich. Finally, we showed that the problem of determining a minimum leader set achieving targeted controllability is NP-hard. Therefore, a heuristic leader selection algorithm was given to compute (minimum) leader sets achieving target control. The algorithm comprises two phases: firstly, it computes a minimum root set of the target set, i.e. a set of vertices from which all target nodes can be reached. Secondly, this minimum root set is greedily extended to a leader set achieving target control. 

Both graph-theoretic conditions for strong targeted controllability provided in this paper are not one-to-one correspondences. Hence, finding a necessary and sufficient topological condition for strong targeted controllability is still an open problem. Furthermore, investigating other system-theoretic concepts like disturbance decoupling and fault detection for the class of distance-information preserving matrices is among the possibilities for future research.

\ifCLASSOPTIONcaptionsoff
 \newpage
\fi

\bibliographystyle{plain}
\bibliography{MyRef}

\begin{thebibliography}{10}

\bibitem{hardness}
A.~Aazami.
\newblock Hardness results and approximation algorithms for some problems on
  graphs.
\newblock {\em PhD thesis, University of Waterloo}, 2008.

\bibitem{BDD}
D.~Bergman, A.~A. Cire, W.~Hoeve, and T.~Yunes.
\newblock Bdd-based heuristics for binary optimization.
\newblock {\em Journal of Heuristics}, 20(2):211--234, 2014.

\bibitem{burgarth2013}
D.~Burgarth, D.~D'Alessandro, L.~Hogben, S.~Severini, and M.~Young.
\newblock Zero forcing, linear and quantum controllability for systems evolving
  on networks.
\newblock {\em IEEE Transactions on Automatic Control}, 58:2349--2354, 2013.

\bibitem{chapman}
A.~Chapman and M.~Mesbahi.
\newblock On strong structural controllability of networked systems: A
  constrained matching approach.
\newblock In {\em American Control Conference (ACC), 2013}, pages 6126--6131,
  June 2013.

\bibitem{egerstedt2012}
M.~Egerstedt, S.~Martini, M.~Cao, K.~Camlibel, and A.~Bicchi.
\newblock Interacting with networks: How does structure relate to
  controllability in single-leader, consensus networks?
\newblock {\em IEEE Control Systems}, 32(4):66--73, Aug 2012.

\bibitem{barabasitarget}
J.~Gao, Y.~Y. Liu, R.~M. D'Souza, and A.~L. Barab\'{a}si.
\newblock {Target control of complex networks}.
\newblock {\em Nature Communications}, 5, November 2014.

\bibitem{glover1976}
K.~Glover and L.~Silverman.
\newblock Characterization of structural controllability.
\newblock {\em IEEE Transactions on Automatic Control}, 21(4):534--537, Aug
  1976.

\bibitem{quantumcontrol}
C.~D. Godsil.
\newblock Control by quantum dynamics on graphs.
\newblock {\em Physical Review A}, 81(5), 2010.

\bibitem{hogben}
L.~Hogben.
\newblock Minimum rank problems.
\newblock {\em Linear Algebra and its Applications}, 432(8):1961--1974, 2010.

\bibitem{lincontrollability}
C.T. Lin.
\newblock Structural controllability.
\newblock {\em IEEE Transactions on Automatic Control}, 19(3):201--208, Jun
  1974.

\bibitem{barabasicontrol}
Y.~Y. Liu, J.~J. Slotine, and A.~L. Barabasi.
\newblock {Controllability of complex networks}.
\newblock {\em Nature}, 473(7346):167--173, May 2011.

\bibitem{strongstructuralcontrollability1979}
H.~Mayeda and T.~Yamada.
\newblock Strong structural controllability.
\newblock {\em SIAM Journal on Control and Optimization}, 17(1):123--138, 1979.

\bibitem{Nima}
N.~Monshizadeh, M.~K. Camlibel, and H.L. Trentelman.
\newblock Strong targeted controllability of dynamical networks.
\newblock In {\em 54th {IEEE} Conference on Decision and Control (CDC), Osaka,
  Japan, December 15-18, 2015}, pages 4782--4787.

\bibitem{zeroforce}
N.~Monshizadeh, S.~Zhang, and M.~K. Camlibel.
\newblock Zero forcing sets and controllability of dynamical systems defined on
  graphs.
\newblock {\em IEEE Transactions on Automatic Control}, 59(9):2562--2567, 2014.

\bibitem{rahmani2009}
A.~Rahmani, M.~Ji, M.~Mesbahi, and M.~Egerstedt.
\newblock Controllability of multi-agent systems from a graph-theoretic
  perspective.
\newblock {\em SIAM Journal on Control and Optimization}, 48(1):162--186, 2009.

\bibitem{faults}
P.~Rapisarda, A.~R.~F. Everts, and M.~K. Camlibel.
\newblock Fault detection and isolation for systems defined over graphs.
\newblock In {\em 54th IEEE Conference on Decision and Control (CDC)}, pages
  3816--3821, Dec 2015.

\bibitem{multiinput}
R.~Shields and J.~Pearson.
\newblock Structural controllability of multiinput linear systems.
\newblock {\em IEEE Transactions on Automatic Control}, 21(2):203--212, Apr
  1976.

\bibitem{tanner2004}
H.~G. Tanner.
\newblock On the controllability of nearest neighbor interconnections.
\newblock In {\em 43rd IEEE Conference on Decision and Control (CDC)},
  volume~3, pages 2467--2472, Dec 2004.

\bibitem{trefois}
M.~Trefois and J.~C. Delvenne.
\newblock Zero forcing number, constrained matchings and strong structural
  controllability.
\newblock {\em Linear Algebra and its Applications}, 484:199--218, 2015.

\bibitem{zhang2011}
S.~Zhang, M.~Cao, and M.~K. Camlibel.
\newblock Upper and lower bounds for controllable subspaces of networks of
  diffusively coupled agents.
\newblock {\em IEEE Transactions on Automatic Control}, 59(3):745--750, 2014.

\end{thebibliography}

\end{document}